\documentstyle[12pt,fleqn]{article}
\begin{document}
\renewcommand{\theequation}{\arabic {section}.\arabic {equation}}
  \newcommand{\bi}{\begin{equation}}
  \newcommand{\ei}{\end{equation}}
   \date{}
\baselineskip 24pt
\title{A Note on  the  Linnik's constant}
\author{{Zaizhao MENG}}
\maketitle
\section{Introduction}
Suppose that $q$ is a sufficiently large positive integer,
$(a,q)=1$, let $P(a,q)$ be the least prime in the arithmetic
progression $\{ n\equiv a(mod\ q)\}$ and $\chi$ be the Dirichlet
character of modulus q. $L(s,\chi)$ is the Dirichlet L-function.
Every positive integer $q$ can be expressed as $q=q_{3}^{3}q_{2},\
q_{2}$ is cube-free. An integer $q$ is called ''has bounded cubic
part'' if the above $ q_{3}$ is bounded
absolutely.\\
In this note, we correct the Lemma of [4] and obtain the following
result.\\
 \ {\bf Theorem }\ If $q$ has bounded cubic part,
 for $(a,q)=1$, we have $$ P(a,q)\ll q^{4.5} .$$
{\bf Notation:} $\varepsilon$  is a sufficiently small positive real
number. $\ell=\log q$, $\chi$ is the Dirichlet character of modulus
$q$ and $\chi_{0}$ is the principal character, $q$ is sufficiently
large positive integer and has bounded cubic part.
\section{Basic lemma and proof of the theorem}
Let $R$ be the rectangle (6.1) of [2]. Let $\rho_{1}$ be a zero of
$\prod\limits_{\chi(mod\ q)}L(s,\chi)$ in $R$ with
$\beta_{1}=\Re(\rho_{1})$ maximal, let $\chi_{1}$ be any character
with $L(\rho_{1},\chi_{1})=0$. Let $\rho_{2}$ be a zero of
$\prod\limits_{\chi\neq\chi_{1}, \bar{\chi_{1}}}L(s,\chi)$ in $R$
with $\beta_{2}=\Re(\rho_{2})$ maximal, let $\chi_{2}$ be any
character with $L(\rho_{2},\chi_{2})=0$,.... We obtain new zeros
$\rho_{1}$,..., and
 $
 \Re\rho_{k+1}\leq\Re\rho_{k}\leq ...\leq\Re\rho_{2}\leq\Re\rho_{1}.
$\\
As in [2], we shall set$
 \rho_{k}=\beta_{k}+i\gamma_{k},\ \beta_{k}=1-\ell^{-1}\lambda_{k},\ \gamma_{k}=
\ell^{-1}\mu_{k}. $\
 The zero $\rho^{\prime}$ of $L(s,\chi)$ is as
that in [2], \ $
 \rho^{\prime}=\beta^{\prime}+i\gamma^{\prime},\ \beta^{\prime}=1-\ell^{-1}\lambda^{\prime},\ \gamma^{\prime}=
\ell^{-1}\mu^{\prime}. $\\
 Using the results in [1] and [2], set
$\phi=\frac{1}{4}$ in the Lemma 3.1 of [2]. The bounds of
$\lambda_{1},\lambda_{2},\lambda_{3}$ and
$\lambda^{\prime}$ in [3] can be used to deduce the theorem. As in [3] and [4], we have
 $\lambda_{3}>\frac{8}{7}-\varepsilon$ except for the case:
$\chi_{1}\chi_{2}\chi_{3}=\chi_{0}, \chi_{1}, \chi_{2},
\chi_{3},\rho_{1}$ all are real.\\
{\bf Lemma. }\ Suppose that $\chi_{1}\chi_{2}\chi_{3}=\chi_{0}$, and
$ \chi_{1}, \chi_{2}, \chi_{3},\rho_{1}$ all are real. If
$\lambda_{1}\leq
0.85$, \ then $\lambda_{3}>\frac{8}{7}-\varepsilon.$\\
{\bf Proof.} As in [4], suppose that $0.6<\lambda_{1}\leq 0.85$. Let
$P_{3}(X)=X+X^{2}+\frac{2}{3}X^{3}$, as
 in [2](on page 315), we can deduce that if $\lambda_{3}\leq b\leq (\sqrt[3]{3}-1)a$, then
\bi
P_{3}(\frac{a+\lambda_{1}}{a})-P_{3}(1)-2P_{3}(\frac{a+\lambda_{1}}{a+b})
+(a+\lambda_{1})(\frac{3}{4}+\varepsilon)\geq 0. \ei We choose
$b=1.15,\ a=b(\sqrt[3]{3}-1)^{-1}$, and will prove that when
$\varepsilon$ is sufficiently small and $q$ is sufficiently large,
(2.1) can not happen,
and then $\lambda_{3}>1.15>\frac{8}{7}$.\\
Let\\
$$
 S(X)=P_{3}(\frac{X}{a})-P_{3}(1)-2P_{3}(\frac{X}{a+b})
+\frac{3}{4}X, $$
 then the differential function
 $$
 S^{\prime}(X)=c_{1}+c_{2}X+c_{3}X^{2},
 $$
where
$$
c_{1}=\frac{3}{4}+\frac{1}{a}-\frac{2}{a+b},
 c_{2}=\frac{2}{a^{2}}-\frac{4}{(a+b)^{2}},
 c_{3}=\frac{2}{a^{3}}-\frac{4}{(a+b)^{3}}.
$$
$S^{\prime}(X)=0$ has not real root, and $S(a+0.6)=-0.269776...,\
S(a+0.85)=-0.0051...$ , when $\varepsilon$ is sufficiently small and
$q$ is sufficiently large, the left hand side of (2.1) is less than
$-0.005$, this contradiction shows that
$\lambda_{3}>1.15$. The lemma follows.\\
Suppose $N(\lambda)$ and $H_{0}$ have the same meaning as in [3].
Taking $\phi=\frac{1}{4}$ in Lemma 12.1 of [2], we add the following
results.\\
$ \lambda_{1}=0.85: N(0.9)=5,\ N(1.0)=6,\ N(1.1)=7,\ N(1.2)=9,\
N(1.3)=12,\ N(1.4)=16,\ N(1.5)=23. $\\
As in [3], we need to show $H_{0}<1$.\\
By the results of [3], it is necessary to add the following case , as in the proof of the Theorem 1.4 of [3]:\\
(a1) $\chi_{1}$ is real, $\rho_{1}$ is real, $\lambda_{1}\in
[0.85,1.15]$.\\
In this case (a1), using the above bounds for $N(\lambda)$, we
deduce that $H_{0}<0.9$. This completes the proof of the theorem.
 
\noindent
{\small School of Mathematical Sciences}\\
{\small Peking University}\\
{\small Beijing, 100871, P.R.China}   \\
{\small E-mail:\ mzzh@math.pku.edu.cn}
\end{document}